\documentclass[preprint,12pt]{elsarticle}

\newtheorem{theorem}{Theorem}[section]
\newtheorem{corollary}{Corollary}[section]
\newtheorem{lemma}{Lemma}[section] 

\usepackage{amssymb}

\begin{document}

\begin{frontmatter}

\title{Large deviation exponential inequalities  for supermartingales}
\author{Xiequan Fan$^*$,\ \ \ \ \ \ Ion Grama,\ \ \ \ \ \ Quansheng Liu}
 \cortext[cor1]{\noindent Corresponding author. \\
\mbox{\ \ \ \ }\textit{E-mail}: fanxiequan@hotmail.com (X. Fan), \ \ \ \ \  ion.grama@univ-ubs.fr (I. Grama),\\
\mbox{\ \ \ \ \ \ \ \ \ \ \  \ \ \ \ }quansheng.liu@univ-ubs.fr (Q. Liu). }
\address{Universit\'{e} de Bretagne-Sud, LMBA, UMR CNRS 6205,
 Campus de Tohannic, 56017 Vannes, France}
\begin{abstract}
Let $(X_{i}, \mathcal{F}_{i})_{i\geq 1}$ be a sequence of  supermartingale differences and let $S_k=\sum_{i=1}^k X_i$. We give an exponential moment condition under which $P( \max_{1\leq k \leq n} S_k \geq n)=O(\exp\{-C_1 n^{\alpha}\}),$ $n\rightarrow \infty, $
where $\alpha \in (0, 1)$ is given and $C_{1}>0$ is a constant. We also show that the power $\alpha$ is optimal under the given condition. In particular, when $\alpha=\frac{1}{3}$, we recover an inequality of  Lesigne and Voln\'{y}.
\end{abstract}
\begin{keyword}Large deviations;  martingales; sub-exponential distributions
\vspace{0.3cm}
\MSC  60F10 \sep 60G42
\end{keyword}

\end{frontmatter}


\section{Introduction}
Let $(X_{i}, \mathcal{F}_{i})_{i\geq 1}$ be a sequence of martingale differences  and let $S_k=\sum_{i=1}^k X_i, k\geq 1$.
Under the Cram\'{e}r condition $\sup_{ i }E e^{|X_{i}|} < \infty$, Lesigne and Voln\'{y} \cite{LV01} proved  that
\begin{eqnarray}
 P(  S_{n}  \geq n ) &=& O (\exp\{ - C_1 n^{\frac{1}{3}}\}), \ \ \ \ \  n \rightarrow \infty,  \label{flv01}
\end{eqnarray}
for some constant $C_1>0$.  Here and throughout the paper, for two functions $f$ and $g$, we write $f(n)=O(g(n))$ if
there exists a constant $C>0$ such that $|f(n)|\leq C|g(n)|$ for all $n\geq1$.
Lesigne and Voln\'{y} \cite{LV01} also showed that the power  $\frac{1}{3}$ in (\ref{flv01}) is optimal in the sense that there exists a sequence of martingale differences $( \widehat{X}_{i}, \mathcal{F}_{i})_{i\geq 1}$
such that $\sup_{i}Ee^{| \widehat{X} _i|}< \infty$ and $P(   \widehat{S} _n \geq n)> \exp\{-C_2 n^{\frac{1}{3}}\}$ for some constant $C_2>0$ and infinitely many $n$'s.
Liu and Watbled \cite{Liu09a} proved that the power $\frac{1}{3}$ in (\ref{flv01}) can be
improved to $1$  under the conditional Cram\'{e}r condition $\sup_{i}E(e^{|X_i|}|\mathcal{F}_{i-1})\leq  C_3$, for some constant $C_3$.
It seems natural to ask under what condition it holds
\begin{eqnarray}
 P(  S_{n}  \geq n ) &=&  O(\exp\{ - C_1 n^{\alpha}\}),\ \ \ \ n\rightarrow \infty,\label{lvie}
\end{eqnarray}
where $\alpha \in (0,1)$ is given and $C_1>0$ is a constant.
In this paper, we give some sufficient conditions in order that (\ref{lvie}) holds for  supermartingales $(S_k,\mathcal{F}_{k})_{k\geq 1}$.

The paper is organized as follows. In Section \ref{sec1}, we present the main results. In Sections \ref{sec2}-\ref{sec4}, we give
the proofs of the main results.

\section{Main Results} \label{sec1}
Our first result is an extension of the bound (\ref{flv01}) of Lesigne and Voln\'{y}.
\begin{theorem}
\label{th1} Let $\alpha \in (0, 1)$. Assume that $(X_{i}, \mathcal{F}_{i})_{i\geq 1}$ is a sequence of
 supermartingale differences  satisfying
 $\sup_{ i }E\exp\{ |X_{i}|^{\frac{2\alpha}{1-\alpha}}\}\leq C_1 $
 for some constant $C_1\in (0, \infty)$. Then, for all $ x>0$,
\begin{eqnarray}
P\left( \max_{1\leq k \leq n}S_k \geq  n x \right) &\leq& C(\alpha, x)  \exp\left\{-\left(\frac{ x  }{ 4 }\right)^{2\alpha} n^\alpha\ \right\} ,\nonumber
\end{eqnarray}
where
\[
C(\alpha,  x)=  2+ 35 C_1  \left( \frac{1}
{ x^{2\alpha} 16^{1-\alpha}}  + \frac{1}{  x^2 }
\left( \frac{3(1-\alpha)}{2\alpha}\right)^{\frac{ 1-\alpha}{\alpha}} \right)
\]
does not depend  on $n$.
In particular, with $x=1$, it holds
\begin{eqnarray}
P\left( \max_{1\leq k \leq n} S_k \geq  n  \right)  &=&O \left( \exp\{ - \frac{1}{16} \, n^{\alpha} \} \right),\ \ \ \ n\rightarrow \infty. \label{fnk}
\end{eqnarray}
 Moreover, the power $\alpha$ in (\ref{fnk}) is optimal even for the class of stationary  martingale differences:
 for each $\alpha \in (0, 1)$, there exists a stationary  sequence of martingale differences $(\widehat{X}_{i}, \mathcal{F}_{i})_{i\geq 1}$ satisfying $E\exp\{ |\widehat{X}_{1}|^{\frac{2\alpha}{1-\alpha}}\}< \infty$ and
\begin{eqnarray}
P\left( \max_{1\leq k \leq n} \widehat{S}_k \geq  n  \right)  &\geq & \exp\{ - 3 n^{\alpha} \},
\end{eqnarray}
for all $n$ large enough.
\end{theorem}

It is clear that when $\alpha=\frac{1}{3}$, the bound (\ref{fnk}) implies the bound (\ref{flv01}) of Lesigne and Voln\'{y}.

In our second result we replace the condition $\sup_{ i }E\exp\{ | X_{i}|^{\frac{2\alpha}{1-\alpha}}\}< \infty $ of Theorem \ref{th1} by the weaker condition $\sup_{i}E \exp\{  ( X_{i}^+)^{\frac{\alpha}{1-\alpha} }\} < \infty$, where  $X_i^+=\max\{X_i, 0\}.$
Denote by $$\langle S\rangle_k= \sum_{i=1}^k E(X_i^2|\mathcal{F}_{i-1})$$  the sum of conditional variances.
\begin{theorem} \label{th2}
Let $\alpha \in (0,1)$. Assume that $(X _i,\mathcal{F}_i)_{i\geq 1}$ is a sequence of
 supermartingale differences  satisfying  $\sup_{i}E \exp\{  ( X_{i}^+)^{\frac{\alpha}{1-\alpha} }\} \leq C_1$ for some constant $ C_1 \in (0, \infty)$. Then,
  for all $x,v>0$,
\begin{eqnarray}
  &&P\left( S_{k}\geq x\ \mbox{and}\ \langle S\rangle_k \leq
 v^2\  \mbox{for some}\ k\in [1,n] \right) \nonumber\\
 &&\quad\quad\quad\quad\quad\quad\quad\quad\quad\quad\quad\quad\quad  \leq \exp\left\{-\frac{x^2}{2(v^2+\frac{1}{3}x^{2-\alpha})}\right\}  + nC_1 \exp\{ -   x^{\alpha } \}.\label{f29}
\end{eqnarray}
\end{theorem}

We note that, for bounded random variables, some inequalities closely related to (\ref{f29}) can be found in  Freedman \cite{FR75},  Dedecker \cite{D01}, Dzhaparidze and van Zanten \cite{Dz01}, Merlev\`{e}de, Peligrad and Rio \cite{R09} and Delyon \cite{D09}.

Adding a hypothesis on $\langle S\rangle_n$ to Theorem \ref{th2},
we can easily obtain the following Bernstein type inequality which is similar to an inequality of Merlev\`{e}de, Peligrad and Rio \cite{R10} for weakly dependent sequences.
\begin{corollary} \label{co1}
Let $\alpha \in (0,1)$.
Assume that $(X _i,\mathcal{F}_i)_{i\geq 1}$ is a sequence of
 supermartingale differences  satisfying
 $E \exp\{ ( X_{i}^+)^{\frac{\alpha}{1-\alpha} }\} \leq C_1$
and
$E\exp\{ (\frac{\langle S\rangle_n}{n})^{\frac{\alpha}{1-\alpha} }\} \leq  C_2$ for some constants $ C_1, C_2 \in (0, \infty)$. Then,
for all  $x> 0$,
\begin{eqnarray}
 P\left(\max_{1 \leq k \leq n} S_{k} \geq nx \right)
&\leq&  \exp\left\{-\frac{  x^{1+\alpha}  }{2\left(1+ \frac{ 1}{3} x \right)} n^{\alpha}  \right\}
  + (nC_1+C_2) \exp\{ - x^{\alpha} n^{\alpha}  \}. \label{fkg}
\end{eqnarray}
In particular, with $x=1$, it holds
\begin{eqnarray}
P\left( \max_{1\leq k \leq n} S_k \geq  n  \right)  &=&O \left( \exp\{ -C  \, n^{\alpha} \} \right),\ \ \ \ n\rightarrow \infty, \label{fkfg}
\end{eqnarray}
where $C  > 0$ is an absolute constant. Moreover, the power $\alpha$ in (\ref{fkfg}) is optimal even for the class of stationary  martingale differences: for each $\alpha \in (0, 1)$, there exists a stationary  sequence of martingale differences $(\widehat{X}_{i}, \mathcal{F}_{i})_{i\geq 1}$ satisfying $E\exp\{ |\widehat{X}_{1}|^{\frac{2\alpha}{1-\alpha}}\}< \infty$ and
\begin{eqnarray}
P\left( \max_{1\leq k \leq n} \widehat{S}_k \geq  n  \right)  &\geq & \exp\{ - 3 n^{\alpha} \},
\end{eqnarray}
for all $n$ large enough.
\end{corollary}

In the i.i.d. case, the conditions of Corollary \ref{co1} can be weakened considerably, see Lanzinger and Stadtm\"{u}ller \cite{LS00}
where it is shown that if $E \exp\{ ( X_{i}^+)^{\alpha  }\} < \infty$ with $ \alpha \in (0, 1)$, then
\begin{eqnarray}
P\left( \max_{1\leq k \leq n} S_k \geq  n  \right)  &=&O \left( \exp\{ -C_\alpha \, n^{\alpha } \} \right),\ \ \ \ n\rightarrow \infty.
\end{eqnarray}

\section{Proof of Theorem \ref{th1}}\label{sec2}
We need the following refined version of the Azuma-Hoeffding inequality.
\begin{lemma} \label{lemma1}
 Assume that $(X_{i}, \mathcal{F}_{i})_{i\geq 1} $ is a sequence of  martingale differences  satisfying $| X_{i}| \leq 1$ for all $i$.
 Then, for all $ x \geq 0$,
\begin{eqnarray}
P\left( \max_{1\leq k \leq n}  S_k  \geq x  \right)&\leq&  \exp\left\{ - \frac{ x^2}{2n}\right\} . \label{hoemax}
\end{eqnarray}
\end{lemma}

A proof can be found in Laib \cite{La99}.

Now, we are ready to prove Theorem \ref{th1}. We start as in Lesigne and Voln\'{y} \cite{LV01} and push a step further
by using the martingale maximal inequality (\ref{hoemax}). We end by giving a simple example to show  that the power $\alpha$ in (\ref{fnk}) is optimal.

Let $(X_{i}, \mathcal{F}_{i})_{i \geq 1}$ be a sequence of supermartingale differences.  Given $u> 0$, define
\begin{eqnarray*}
X'_{i} &=& X_{i}\mathbf{1}_{\{|X_{i}|\leq u\}}-E(X_{i}\mathbf{1}_{\{|X_{i}|\leq u\}}|\mathcal{F}_{i-1}),\\
X''_{i} &=& X_{i}\mathbf{1}_{\{|X_{i}|> u\}}-E(X_{i}\mathbf{1}_{\{|X_{i}|> u\}}|\mathcal{F}_{i-1}), \\
S'_k &=&\sum_{i=1}^k X'_{i},  \ \ \ \ \ S''_k =\sum_{i=1}^k X''_{i},  \ \ \ \ S'''_k = \sum_{i=1}^k  E(X_{i}|\mathcal{F}_{i-1}).
\end{eqnarray*}
Then $(X'_{i}, \mathcal{F}_{i})_{i\geq 1}$ and $(X''_{i}, \mathcal{F}_{i})_{i\geq 1}$ are two martingale difference sequences and
$S_k= S'_k + S''_k+S'''_k$.
Let $t \in (0,1).$ Since $S'''_k\leq 0$, for any $x>0$,
\begin{eqnarray}
P\left( \max_{1\leq k \leq n}S_k \geq x \right) &\leq& P\left( \max_{1\leq k \leq n}S'_k+S'''_k \geq x t \right) + P\left( \max_{1\leq k \leq n}S''_k\geq x (1-t) \right) \nonumber\\
&\leq& P\left( \max_{1\leq k \leq n}S'_k \geq x t \right) + P\left( \max_{1\leq k \leq n}S''_k\geq x (1-t) \right). \label{fsum}
\end{eqnarray}
Using Lemma \ref{lemma1} and $|X'_{i}|\leq 2u$, we have
\begin{eqnarray}
P\left( \max_{1\leq k \leq n} S'_k \geq x t \right)&\leq&  \exp\left\{-\frac{ x^2t^2}{8 n u^2} \right\} . \label{f7i}
\end{eqnarray}
Let $F_{i}(x)=P( |X_{i}| \geq x), x\geq 0.$ Since $E\exp\{|X_{i}|^{\frac{2\alpha}{1-\alpha}}\}\leq C_1$, we obtain, for all $x\geq0$,
\[
F_{i}(x)\leq \exp\{ -x^{\frac{2\alpha}{1-\alpha}}\}E\exp\{|X_{i}|^{\frac{2\alpha}{1-\alpha}}\} \leq C_1 \exp\{ -x^{\frac{2\alpha}{1-\alpha}}\}.
\]
Using the martingale maximal inequality p. 14 in \cite{HH80}, we get
\begin{eqnarray}
P\left( \max_{1\leq k \leq n} S''_k \geq x (1-t) \right) &\leq& \frac{1}{x^2(1-t)^2}\sum_{i=1}^n E  X''_{i}\mbox{}^{2}. \label{fk}
\end{eqnarray}
It is easy to see that
\begin{eqnarray}
 EX''_i\mbox{}^2 &=& -  \int_{u}^\infty  t^2dF_{i}(t)\nonumber\\
  &= &  u^2F_{i}(u) + \int_{u}^\infty 2 tF_{i}(t) dt  \nonumber\\
&\leq & C_1 u^2\exp\{ -u^{\frac{2\alpha}{1-\alpha}}\}+2C_1  \int_{u}^\infty t \exp\{ -t^{\frac{2\alpha}{1-\alpha}}\} dt. \label{fkl}
\end{eqnarray}
Notice that the function $g(t)=t^3 \exp\{ -t^{\frac{2\alpha}{1-\alpha}}\}$ is decreasing in $[\beta, +\infty)$ and is increasing in $[0, \beta]$, where $\beta =\left(\frac{3(1-\alpha)}{2\alpha} \right)^{\frac{1-\alpha}{2\alpha}}$.
If $0<  u < \beta$, we have
\begin{eqnarray}
 \int_{u}^\infty t \exp\{ -t^{\frac{2\alpha}{1-\alpha}}\} dt &\leq &  \int_{u}^\beta t \exp\{ -t^{\frac{2\alpha}{1-\alpha}}\} dt   + \int_{\beta}^\infty t^{-2} t^3 \exp\{ -t^{\frac{2\alpha}{1-\alpha}}\} dt \nonumber \\
 &\leq&  \int_{u}^\beta t \exp\{ -u^{\frac{2\alpha}{1-\alpha}}\} dt  + \int_{\beta}^\infty t^{-2} \beta^3 \exp\{ -\beta^{\frac{2\alpha}{1-\alpha}}\} dt \nonumber \\
 &\leq& \frac{3}{2}\beta^2\exp\{ -u^{\frac{2\alpha}{1-\alpha}}\}. \label{ghds}
\end{eqnarray}
If $ \beta \leq  u$, we have
\begin{eqnarray}
 \int_{u}^\infty t \exp\{ -t^{\frac{2\alpha}{1-\alpha}}\} dt &=&  \int_{u}^\infty t^{-2} t^3 \exp\{ -t^{\frac{2\alpha}{1-\alpha}}\} dt \nonumber \\
 &\leq&  \int_{u}^\infty t^{-2} u^3 \exp\{ -u^{\frac{2\alpha}{1-\alpha}}\} dt \nonumber \\
 &=& u^2\exp\{ -u^{\frac{2\alpha}{1-\alpha}}\}. \label{xvcb}
\end{eqnarray}
Returning to (\ref{fkl}), by (\ref{ghds}) and (\ref{xvcb}), we get
\begin{eqnarray}
 EX''_i\mbox{}^2 &\leq &3 C_1( u^2 +  \beta^2 ) \exp\{ -u^{\frac{2\alpha}{1-\alpha}}\}.
\end{eqnarray}
From (\ref{fk}), it follows that
\begin{eqnarray}
P\left( \max_{1\leq k \leq n} S''_k \geq x (1-t) \right)&\leq& \frac{3 n C_1  }{x^2(1-t)^2} (u^2+ \beta^2)\exp\{ -u^{\frac{2\alpha}{1-\alpha}}\}. \label{fl}
\end{eqnarray}
Combining (\ref{fsum}), (\ref{f7i}) and (\ref{fl}), we obtain
\begin{eqnarray}
P\left( \max_{1\leq k \leq n} S_k \geq  x \right) &\leq&   2\exp\left\{-\frac{ x^2t^2}{8 n u^2} \right\} +\frac{3 n C_1  }{ (1-t)^2} \left( \frac{u^2}{x^2}+ \frac{\beta^2}{x^2} \right) \exp \{ -u^{\frac{2\alpha}{1-\alpha}}  \}.\nonumber
\end{eqnarray}
Taking $t=\frac{1}{\sqrt{2}}$ and $u=  \left(\frac{x }{4\sqrt{n}} \right)^{1-\alpha},$ we get, for all $x>0$,
\begin{eqnarray*}
P\left( \max_{1\leq k \leq n} S_k \geq  x \right) &\leq& C_n(\alpha,  x)\exp\left\{-\left(\frac{ x^2 }{ 16 n }\right)^\alpha\ \right\}, \nonumber
\end{eqnarray*}
where
\[
C_n(\alpha,   x) =  2+ 35 n C_1  \left( \frac{1}{x^{2\alpha}(16n)^{1-\alpha}}  + \frac{\beta^2}{x^2} \right) .
\]
Hence, for all $x>0$,
\begin{eqnarray*}
P\left( \max_{1\leq k \leq n} S_k \geq n x \right) &\leq&C(\alpha, x)  \exp\left\{-\left(\frac{ x  }{ 4 }\right)^{2\alpha} n^\alpha\ \right\}, \nonumber
\end{eqnarray*}
where
\[
C(\alpha,  x)=  2+ 35 C_1  \left( \frac{1}
{ x^{2\alpha} 16^{1-\alpha}}  + \frac{1}{  x^2 }
\left( \frac{3(1-\alpha)}{2\alpha}\right)^{\frac{ 1-\alpha}{\alpha}} \right).
\]
This completes the first assertion of Theorem \ref{th1}.

Next, we prove that the power $\alpha$ in (\ref{fnk}) is optimal. We take a positive random variable $X$
such that, for all $x>1$,
\begin{eqnarray}
P\left( X \geq  x \right)= \frac{2 e}{ 1+x^{\frac{1+\alpha}{1-\alpha}} } \exp\left\{ -x^{\frac{2\alpha}{1-\alpha}} \right\}. \label{fghdsg}
\end{eqnarray}
It is easy to verify that
\begin{eqnarray*}
E\exp\{ |X |^{\frac{2\alpha}{1-\alpha}}\} =  -\int_1^{\infty} \exp\{ t^{\frac{2\alpha}{1-\alpha}}\} d P\left( X \geq  t \right) =
e+ \frac{4 e\, \alpha}{1-\alpha} \int_1^{\infty} \frac{t^{\frac{3\alpha-1}{1-\alpha}}}{1 + t^{\frac{1+\alpha}{1-\alpha}}} dt < \infty.
\end{eqnarray*}
Assume that $(\xi_i)_{i\geq 1}$ are Rademacher random variables independent of $X$, i.e. $P(\xi_i=1)=P(\xi_i=-1)=\frac12$. Set $\widehat{X}_i=X\xi_i$ and $\mathcal{F}_i=\sigma(X, (\xi_k)_{k=1,...,i})$. Then, $(\widehat{X}_i, \mathcal{F}_i)_{i\geq 1}$ is a stationary sequence of martingale differences satisfying
$ \sup_{i} E\exp\{ |\widehat{X}_i |^{\frac{2\alpha}{1-\alpha}}\}=E\exp\{ |X |^{\frac{2\alpha}{1-\alpha}}\}  < \infty$. For $\beta \in (0,1)$,  it is easy to see that
\begin{eqnarray*}
P\left(  \max_{1\leq k \leq n} \widehat{S}_i \geq n  \right)  \geq  P\left(   \widehat{S}_n \geq n  \right)
 \geq   P\left(  \sum_{i=1}^n \xi_i \geq n^{\beta}  \right)  P \left(   X \geq n^{1-\beta}  \right).
\end{eqnarray*}
Since, for $n$ large enough,
\[
P\left(  \sum_{i=1}^n \xi_i \geq n^\beta \right) \geq   \exp \left\{ -n^{2\beta-1}  \right\},
\]
(cf. Corollary 3.5 in Lesigne and Voln\'{y} \cite{LV01}),
we get, for $n$ large enough,
\begin{eqnarray}
P\left(  \max_{1\leq k \leq n} \widehat{S}_i \geq n  \right) &\geq& \frac{2 e}{ 1+(n^{1-\beta})^{\frac{1+\alpha}{1-\alpha}} } \exp\left\{-n^{2\beta-1}  -(n^{1-\beta})^{\frac{2\alpha}{1-\alpha}} \right\}.\  \ \
\end{eqnarray}
Setting $2\beta-1=\alpha$, we obtain, for $n$ large enough,
\begin{eqnarray*}
P\left(  \max_{1\leq k \leq n} \widehat{S}_i \geq n  \right)  \geq   \frac{2 e}{ 1+ n^{\frac{1+\alpha}2 } }  \exp \left\{ -2 n^{\alpha} \right\}
 \geq   \exp \left\{ -3 n^{\alpha} \right\},
\end{eqnarray*}
which proves that the power $\alpha$ in (\ref{fnk}) is optimal.

\section{Proof of Theorem \ref{th2}}\label{sec3}

To prove Theorem \ref{th2}, we need the following inequality whose proof can be found in Fan, Grama and Liu \cite{F12}.
\begin{lemma}
\label{lem1} Assume that $(X _i,\mathcal{F}_i)_{i\geq 1}$ are supermartingale differences satisfying $X_{i} \leq 1$ for all $i$. Then,
for all $ x, v > 0$,
\begin{eqnarray}
  P \left( S_k \geq x\ \mbox{and}\ \langle S\rangle_{k}\leq v^2\ \mbox{for some}\ k \in [1, n] \right)  \leq   \exp\left\{- \frac{x^2}{2(v^2+ \frac{1}{3}x)} \right\}.\label{fgl1}
\end{eqnarray}
\end{lemma}

Assume that $(X _i,\mathcal{F}_i)_{i\geq 1}$ are supermartingale differences.
Given $u> 0$, set
\begin{eqnarray*}
X'_{i}  =   X_{i}\mathbf{1}_{\{ X_{i} \leq u\}},\ \ X''_{i} = X_{i}\mathbf{1}_{\{ X_{i} > u\}},\ \ S'_k  = \sum_{i=1}^k X'_{i} \ \ \mbox{and} \ \ S''_k = \sum_{i=1}^k X''_{i}.
\end{eqnarray*}
Then, $(X'_{i}, \mathcal{F}_{i})_{i\geq 1}$ is also a sequence of supermartingale differences and
$S_k= S'_k + S''_k$. Since $\langle S'\rangle_k \leq \langle S\rangle_k$, we deduce, for any $x, u, v>0$,
\begin{eqnarray}
&&P\left( S_k \geq x\ \mbox{and}\ \langle S\rangle_k\leq v^2 \ \mbox{for some}\ k \in [1, n] \right)\nonumber \\
&\leq&P \left( S_k' \geq x\ \mbox{and}\ \langle S\rangle_k\leq v^2 \ \mbox{for some}\ k\in [1, n]  \right) \nonumber \\
& & + P\left( S_k'' \geq 0\ \mbox{and}\ \langle S\rangle_k\leq v^2 \ \mbox{for some}\ k \in [1, n] \right) \nonumber \\
&\leq&P\left( S_k' \geq x\ \mbox{and}\ \langle S'\rangle_k\leq v^2 \ \mbox{for some}\ k \in [1, n] \right)   + P \left( \max_{1\leq k \leq n} S_k'' \geq 0  \right).  \label{kdjn}
\end{eqnarray}
 Applying Lemma \ref{lem1} to the supermartingale differences $\left( X'_i/u,\mathcal{F}_i \right)_{i\geq 1}$, we have
\begin{eqnarray}
 P( S_k' \geq x\ \mbox{and}\ \langle S'\rangle_k\leq v^2 \ \mbox{for some}\ k \in [1, n] )
 \leq  \exp\left\{-\frac{x^2}{2(v^2+\frac{1}{3}xu)}\right\}. \label{sgfd}
\end{eqnarray}
Using the exponential Markov's inequality and the condition $E \exp\{ (X_{i}^+)^{\frac{\alpha}{1-\alpha} }\}$ $\leq C_1$, we get
\begin{eqnarray}
P\left( \max_{1\leq k \leq n} S_k'' \geq 0  \right) &\leq& \sum_{i=1}^n P(X_{i}>u) \nonumber\\
 &\leq& \sum_{i=1}^n E \exp\{   (X_{i}^+)^{\frac{\alpha}{1-\alpha} }-u^{\frac{\alpha}{1-\alpha}} \} \nonumber\\
&\leq& nC_1 \exp\{ -   u^{\frac{\alpha}{1-\alpha} } \}. \label{f45}
\end{eqnarray}
Combining the inequalities (\ref{kdjn}), (\ref{sgfd}) and (\ref{f45}) together, we obtain, for all $x, u, v>0$,
\begin{eqnarray}
&& P( S_k \geq x\ \mbox{and}\ \langle S\rangle_k\leq v^2 \ \mbox{for some}\ k \in [1,n] )\nonumber\\
 &&\quad\quad\quad\quad\quad\quad\quad\quad\quad \quad\quad    \leq  \exp\left\{-\frac{x^2}{2(v^2+\frac{1}{3}xu)}\right\}  + n C_1 \exp\{ -   u^{\frac{\alpha}{1-\alpha} } \}.
\end{eqnarray}
Taking $u= x^{1-\alpha}$, we get, for all $x,v>0$,
\begin{eqnarray}
&& P( S_k \geq x\ \mbox{and}\ \langle S\rangle_k\leq v^2 \ \mbox{for some}\ k \in [1,n] )\nonumber\\
  &&\quad\quad\quad\quad\quad\quad\quad\quad\quad \quad\quad \leq  \exp\left\{-\frac{x^2}{2(v^2+\frac{1}{3}x^{2-\alpha})}\right\}  + nC_1 \exp\{ -  x^{\alpha } \}.
\end{eqnarray}
This completes the proof of Theorem \ref{th2}.

\section{Proof of Corollary \ref{co1}.}\label{sec4}
To prove Corollary \ref{co1} we make use of Theorem \ref{th2}.  It is easy to see that
\begin{eqnarray}
P\left(\max_{1\leq k \leq n} S_{k}\geq nx \right) &\leq& P\left(\max_{1\leq k \leq n} S_{k}\geq nx,  \langle S\rangle_n \leq nv^2\right) \nonumber\\  &&+ P\left(\max_{1\leq k \leq n} S_{k}\geq nx,  \langle S\rangle_n > nv^2\right) \nonumber\\
&\leq& P( S_k \geq nx\ \mbox{and}\ \langle S\rangle_k\leq nv^2 \ \mbox{for some}\ k \in [1,n] ) \nonumber\\
 && +  P\left( \langle S\rangle_n > nv^2\right).
\end{eqnarray}
By Theorem \ref{th2}, it follows that, for all $x, v>0$,
\begin{eqnarray}
 P\left(\max_{1\leq k \leq n} S_{k}\geq nx \right)
&\leq&  \exp\left\{-\frac{  x^2 }{2\left( n^{\alpha-1}v^2+ \frac{1}{3} x^{2-\alpha} \right)}n^{\alpha} \right\} \nonumber\\
& &  + n C_1 \exp\left\{ -  x^{\alpha} n^{\alpha} \right\}+ P(\langle S\rangle_n > nv^2 ), \nonumber
\end{eqnarray}
Using the exponential Markov's inequality and the condition $E\exp\{ (\frac{\langle S\rangle_n}{n})^{\frac{\alpha}{1-\alpha} }\} \leq  C_2$, we get, for all $ v>0$,
\begin{eqnarray}
P\left(\langle S\rangle_n > nv^2 \right)&\leq& E \exp\left\{ \left((\frac{\langle S\rangle_n}{n})^{\frac{\alpha}{1-\alpha}}-v^{2\frac{\alpha}{1-\alpha}} \right)\right\}
\leq C_2 \exp\{ -  v^{2\frac{\alpha}{1-\alpha} } \}. \nonumber
\end{eqnarray}
Taking $v =(nx )^{\frac{1-\alpha}{2}}$, we obtain, for all $x>0$,
\begin{eqnarray}
 P\left(\max_{1 \leq k \leq n} X_{k}\geq nx \right)
&\leq&  \exp\left\{-\frac{  x^{1+\alpha}  }{2\left( 1+ \frac{1}{3} x \right)} n^{\alpha}  \right\} + (nC_1+C_2) \exp\{ -  x^{\alpha} n^{\alpha}  \},\nonumber
\end{eqnarray}
which gives  inequality (\ref{fkg}).

Next, we prove that the power $\alpha$ in (\ref{fkfg}) is optimal even for the class of stationary  martingale differences. Let $X$ be the positive random variable defined in (\ref{fghdsg}).
Let $\widehat{X}_i=X\xi_i$ and $\mathcal{F}_i=\sigma(X, (\xi_k)_{k=1,...,i})$, where $(\xi_i)_{i\geq 1}$ are Rademacher random variables independent of $X$. Note that $\frac{\langle \widehat{S}\rangle_n}{n}=X^2$ satisfies the condition $$E\exp\{ (\frac{\langle \widehat{S}\rangle_n}{n})^{\frac{\alpha}{1-\alpha} }\} = E\exp\{|X|^{\frac{2\alpha}{1-\alpha} }\} < \infty.$$
Using the same argument as in the proof of Theorem \ref{th1}, we obtain, for $n$ large enough,
\begin{eqnarray*}
P\left(  \max_{1\leq k \leq n} \widehat{S}_k \geq n  \right) &\geq&  \exp \left\{ -3 n^{\alpha} \right\},
\end{eqnarray*}
which shows that the power $\alpha$ in (\ref{fkfg}) is optimal.


\begin{thebibliography}{99}{\footnotesize
\bibitem{D01}  Dedecker, J., 2001. Exponential inequalities and functional central limit theorems for random fields.
\emph{ESAIM Probab. Statit.} \textbf{5}, 77--104.
\bibitem{D09} Delyon, B., 2009. Exponential inequalities for sums of weakly dependent variables.
\emph{Electron. J. Probab.} \textbf{14},  752--779.
\bibitem{Dz01} Dzhaparidze, K. and van Zanten, J. H.,  2001.  On Bernstein-type inequalities for martingales.
\textit{Stochastic Process. Appl.} \textbf{93}, 109--117.
\bibitem{F12} Fan, X., Grama, I. and Liu, Q.,  2012.
 Hoeffding's inequality for supermartingales. \emph{Stochastic Process. Appl.} \textbf{122}, 3545--3559.
\bibitem{FR75}Freedman, D. A., 1975. On tail probabilities
for martingales. \emph{Ann. Probab.} \textbf{3}, 100--118.
\bibitem{HH80} Hall, P. and Heyde, C. C.,  1980.
\emph{Martingale Limit Theory and Its Application}, Academic Press, 81--96.
\bibitem{La99} Laib, N., 1999. Exponential-type inequalities for martingale difference sequences. Application
to nonparametric regression estimation. \emph{Commun. Statist.-Theory. Methods}, \textbf{28}, 1565--1576.
\bibitem{LS00} Lanzinger, N. and Stadtm\"{u}ller, U., 2000.  Maxima of increments of partial sums for certain subexponential distributions. \emph{Stochastic Process. Appl.}, \textbf{86}, 307--322.
\bibitem{LV01} Lesigne, E. and Voln\'{y}, D., 2001. Large
deviations for martingales. \emph{Stochastic Process. Appl.} \textbf{96},
143--159.
\bibitem{Liu09a} Liu, Q. and Watbled, F., 2009. Exponential
ineqalities for martingales and asymptotic properties of the free energy of
directed polymers in a random environment. \emph{Stochastic Process. Appl.}
\textbf{119}, 3101--3132.
\bibitem{R10}  Merlev\`{e}de, F., Peligrad, M. and Rio, E., 2011.
 A Bernstein type inequality and moderate deviations for weakly
dependent sequences. \emph{Probab. Theory Relat. Fields} \textbf{151}, 435--474.
\bibitem{R09}  Merlev\`{e}de, F., Peligrad, M. and Rio, E., 2009.
Bernstein inequality and moderate deviations under strong mixing conditions. \emph{IMS Collections. High Dimensional Probability} \textbf{5}, 273--292.}
\end{thebibliography}
\end{document}